&\usepackage[T2A]{fontenc}

\documentclass[a4paper,11pt]{amsart}
\usepackage[cp1251]{inputenc}
\usepackage[english]{babel}
\usepackage{amsmath,amssymb,euscript,amsthm,amsfonts}

\newtheorem{lem}{Lemma}
\newtheorem{cor}{Corollary}
\newtheorem{prop}{Proposition}
\newtheorem{theor}{Theorem}
\begin{document}

\title{Minkowski sum of a Voronoi parallelotope and a segment}

\author{Robert Erdahl}
\address{Robert Erdahl, Queen's University, Kingston, Canada}
\email{robert.erdahl@queensu.ca}

\author{Viacheslav Grishukhin}
\address{Viacheslav Grishukhin, CEMI Russian Academy of Sciences,
Nakhimovskii prosp.47 117418 Moscow, Russia}
\email{vgrishukhin@mail.ru}

\date{}
\maketitle
\begin{abstract}
By a {\em Voronoi parallelotope} $P(a)$ we mean a parallelotope determined by a non-negative quadratic form $a$. It was studied by Voronoi in his famous memoir. For a set of vectors $\mathcal P$, we call its {\em dual} a set of vectors ${\mathcal P}^*$ such that $\langle p,q\rangle\in\{0,\pm 1\}$ for all $p\in{\mathcal P}$ and $q\in{\mathcal P}^*$. We prove that Minkowski sum of a Voronoi parallelotope $P(a)$ and a segment is a Voronoi parallelotope $P(a+a_e)$ if and only if this segment is parallel to a vector $e$ of the dual of the set of normal vectors of all facets of $P(a)$, where $a_e(p)=b\langle e,p\rangle^2$ is a quadratic form of rank 1 related to the segment. \end{abstract}

\section{Introduction}

Consider a centrally symmetric $d$-polytope $P(a)$ given by the following system of inequalities
\begin{equation}
\label{Pa}
P(a)=\{x\in{\mathbb R}^d:\langle p,x\rangle\le a(p)\mbox{  for all   } p\in{\mathcal P} \},
\end{equation}
where $\langle p,x\rangle$ is a scalar product of vectors $p,x\in{\mathbb R}^d$. Here ${\mathcal P}\subset{\mathbb R}^d$ is a symmetric set of vectors containing {\em normal} vectors of all facets of  $P(a)$, where {\em symmetric} means that if $p\in{\mathcal P}$ then $-p\in\mathcal P$, too. The function $a:{\mathcal P}\to{\mathbb R}$ is an arbitrary function.

The above $d$-dimensional polytope $P(a)$ is called a {\em Voronoi parallelotope} if the following conditions hold:

(i) the function $a(p)=\langle p,Ap\rangle$ is a non-negative quadratic form;

(ii) the set $\mathcal P$ contains the set ${\mathcal P}_s(a)\subset L$ of normal vectors of all facets of the polytope $P(a)$;

(iii) the set ${\mathcal P}_s(a)$ generates integrally a $d$-dimensional lattice $L$.

Recall that a {\em parallelotope} is a polytope whose parallel translations fill its space without interstices (gaps) and intersections by inner points. Voronoi proved in \cite{Vo} that if the above conditions (i), (ii) and (iii) hold, then $P(a)$ is a parallelotope. Besides, the parallelotope $P(a)$ is a Dirichlet-Voronoi cell of the lattice $2AL$ with respect to the metric form $a$.

One can prove that the set $\mathcal P$ can be enlarged up to a set ${\mathcal P}(a)\subset L$ of minimal (with respect to the form $a$) vectors of each parity class of $L$. Moreover, the set $\mathcal P$ may be the whole lattice $L$.

Let $q\in{\mathbb R}^d$ be a vector and $\alpha\in\mathbb R$ be a number. Define the following affine hyperplanes
\begin{equation}
\label{Hpa}
H(q,\alpha)=\{x\in{\mathbb R}^d:\langle q,x\rangle=\alpha\}\mbox{  and  }
H_p(a)=H(p,a(p)).
\end{equation}

Note that only for $p\in{\mathcal P}(a)$ the hyperplane $H_p(a)$ supports the Voronoi parallelotope $P(a)$ at a face $F(p)$. This face $F(p)$ is called {\em contact face} of $P(a)$. Hence, we call vectors $p\in{\mathcal P}(a)$ by {\em contact vectors}. Dolbilin call in \cite{Do} contact faces by {\em standard} faces.

For each $p\in{\mathcal P}(a)$, the vector $2Ap$ is called {\em commensurate} (with the parallelotope $P(a)$). The commensurate vector $2Ap$ connects the center of the parallelotope $P(a)$ with the center of a parallelotope that is adjacent to $P(a)$ by the contact face $F(p)$. Commensurate vectors generate the lattice $2AL$.

Recall that the dual of a lattice $L$ is
\[L^*=\{q\in{\mathcal R}^d:\langle q,p\rangle\in{\mathbb Z}\mbox{  for all }p\in L\}. \]
Since the set ${\mathcal P}_s(a)$ generate the above lattice $L$, we can change $L$ by ${\mathcal P}_s(a)$ in the above definition of $L^*$. Define the following important subset ${\mathcal P}_s^*(a)\subset L^*$ as follows
\begin{equation}
\label{P*}
{\mathcal P}_s^*(a)=\{e\in{\mathbb R}^d:\langle e,p\rangle\in\{0,\pm 1\}\mbox{  for all  }p\in{\mathcal P}_s(a)\}.
\end{equation}
We call this set {\em dual} of ${\mathcal P}_s(a)$.
Each vector $e\in{\mathcal P}_s^*(a)$ determines a partition of the lattice $L$ onto $(d-1)$-dimensional layers.
\begin{lem}
\label{lay}
Let $e\in{\mathcal P}_s^*(a)$. Then
\[L=\cup_{z\in{\mathbb Z}}L_e(z), \]
where
\begin{equation}
\label{Lez}
L_e(z)=L\cap H(e,z)
\end{equation}
is a $(d-1)$-dimensional layer of $L$, and the hyperplane $H(e,z)$ is defined in (\ref{Hpa}).
\end{lem}
{\bf Proof}. Since the set of normal vectors ${\mathcal P}_s(a)$ generates the lattice $L$, for any $v\in L$, we have $v=\sum_{p\in{\mathcal P}_s(a)} z_p p$, where $z_p\in{\mathbb Z}$. Hence
\[\langle e,v\rangle=\sum_{p\in{\mathcal P}_s(a)}z_p\langle e,p\rangle=z\in{\mathbb Z},\]
i.e. $v\in L_e(z)$. Since $0\in L_e(0)\subset L$, $L_e(z)=v+L_e(0)$. Since $L$ is a $d$-dimensional lattice, each layer has dimension $d-1$. \hfill $\Box$

\vspace{2mm}

Let $e\in{\mathbb R}^d$ be a vector and $l(e)$ be a line spanned by $e$. Let
\begin{equation}
\label{ze}
z(e)=\{x\in{\mathbb R}^d:x=\lambda e, \mbox{  }-1\le\lambda\le 1\}.
\end{equation}
be a segment of the line $l(e)$ symmetric with respect origin 0. For the vector $e$ and a number $b>0$, define the following quadratic form of rank 1
\begin{equation}
\label{ae}
a_e(p)=b\langle p,e\rangle^2,
\end{equation}
We prove below, that $bz(e)=P(a_e)$ if $\langle p,e\rangle\in\{0,\pm 1\}$ for all $p\in{\mathcal P}$.

A parallelotope $P$ is called {\em reducible} if $P=P_1\oplus P_2$, where $\oplus$ denotes direct sum. Otherwise, $P$ is called {\em irreducible}.

In this paper we prove the following
\begin{theor}
\label{thm}
Let $P(a)$ be an irreducible Voronoi parallelotope, defined in (\ref{Pa}), where ${\mathcal P}\supseteq{\mathcal P}(a)\supseteq{\mathcal P}_s(a)$. Let $e\in{\mathbb R}^d$ be a vector. Then one can choose a length of the vector $e$ such that the following assertions are equivalent:

(i) Minkowski sum $P(a)+bz(e)$ is a Voronoi parallelotope for any $b\ge 0$, and
\[P(a)+bz(e)=P(a)+P(a_e)=P(a+a_e);\]

(ii) $e\in{\mathcal P}_s^*(a)$;
\end{theor}

The implication $(ii)\Rightarrow(i)$ was proved in \cite{Gr}.

Magazinov proved in \cite{Ma} that (in our terms) $P(a)+bz(e)$ is a Voronoi parallelotope if this sum is a parallelotope. It seems to us that our proof is simpler.

\section{Minkowski sum of polytopes}
For a fixed set $\mathcal P$ of normal vectors, the polytopes $P(a)$ defined in (\ref{Pa}) have the following simple property
\begin{lem}
\label{a12}
For any functions $a_1(p)$ and $a_2(p)$, the following inclusion holds.
\[P(a_1)+P(a_2)\subseteq P(a_1+a_2). \]
\end{lem}
{\bf Proof}. For $k\in\{1,2\}$, let $x_k\in P(a_k)$. Then $\langle p,x_k\rangle\le a_k(p)$ for all $p\in\mathcal P$. This implies that $\langle p,(x_1+x_1)\rangle\le a_1(p)+a_2(p)$ for all $p\in\mathcal P$, i.e., $x_1+x_2\in P(a_1+a_2)$. Hence $P(a_1)+P(a_2)\subseteq P(a_1+a_2)$. \hfill $\Box$

\vspace{2mm}
Let a set $\mathcal P$ of normal vectors be fixed. It is a problem to find conditions when the equality $P(a_1)+P(a_2)=P(a_1+a_2)$ holds. We show below that this equality holds when $P(a_2)$ is a segment $bz(e)$ and $a_2=f_e$, where the function $f_e(p)$ is defined below in (\ref{fe}).

For $i=1,2$, let $a_i$ be a non-negative quadratic form and $P(a_i)$ be the corresponding Voronoi parallelotope described in (\ref{Pa}).
It is also a problem to find conditions when the sum $P(a_1)+P(a_2)$ is a parallelotope, and, in particular, it is a Voronoi parallelotope.

It is shown in \cite{RBo} that the equality $P(a_1)+P(a_2)=P(a_1+a_2)$ holds if $a_1$ and $a_2$ belong to closure of an L-type domain. We show below that this equality holds if $a_2(p)=a_e(p)$, where the quadratic form $a_e=b\langle p,e\rangle^2$ of rank 1 relates to a segment $bz(e)$, and then the sum $P(a_1)+P(a_e)$ is a parallelotope.

\section{Segments}

Let $e,p\in{\mathbb R}^d$ be some vectors. Consider the affine hyperplane $H_p(f_e)$ defined in (\ref{Hpa}), where
\begin{equation}
\label{fe}
f_e(p)=b\frac{\langle p,e\rangle^2}{|\langle p,e\rangle|}.
\end{equation}
Here $b$ is a non-negative weight of the segment $z(e)$ defined in (\ref{ze}). It is natural to suppose that $f_e(p)=0$ if $\langle p,e\rangle=0$.
\begin{lem}
\label{l1}
For any vector $p\in{\mathbb R}^d$, the hyperplane $H_p(f_e)$ supports the segment $bz(e)$.
\end{lem}
{\bf Proof}. Note that end-vertices of the segment $bz(e)$ are points $\pm be$. If $\langle p,e\rangle>0$, then the end-vertex $be$ lies on $H_p(f_e)$. If $\langle p,e\rangle<0$, then the end-vertex $-be$ lies on $H_p(f_e)$. If $\langle p,e\rangle=0$, then the whole segment $bz(e)$ lies on $H_p(f_e)$. \hfill $\Box$

\vspace{2mm}
Lemma~\ref{l1} implies the following fact.
\begin{lem}
\label{l2}
Let $\mathcal P\subset{\mathbb R}^d$ be a set of vectors such that scalar products $\langle p,e\rangle$ have all the three signs $+,-$ and 0. Let $P(a)$ be given by (\ref{Pa}). Then
\[bz(e)=P(f_e),\]
where the function $f_e(p)$ is defined in (\ref{fe}). \hfill $\Box$
\end{lem}


\section{Minkowski sum of a polytope with a segment}
At first, we consider the Minkowski sum $P(a)+z(e)$ of an {\em arbitrary polytope} $P(a)$ defined in (\ref{Pa}) and the segment $z(e)$ defined in (\ref{ze}).  \'{A}.Horv\'{a}th call in \cite {Ho} the sum $P+z(e)$ by an {\em extension} $P^e$ of $P$. So, we consider a polytope $P=P(a)$ described by the inequalities in (\ref{Pa}), where $a=a(p)$ is an arbitrary function defined on a symmetric set $\mathcal P$. Recall that we call a face $F$ contact and denote it by $F(p)$ if $F=P(a)\cap H_p(a)$. The vector $p$ is called contact vector of the face $F=F(p)$.

For a face $F$ of a polytope $P=P(a)$, let $l_F(e)$ be a parallel shift of the line $l(e)$ such that $l_F(e)\cap F\not=\emptyset$. Call the face $F$ {\em transversal} to $e$ if $l_F(e)\cap F$ is a point. Otherwise, call the face $F$ {\em parallel} to $e$ and denote this fact as $F\parallel e$.

We say that a face $F$ belongs to a {\em shadow boundary} of $P$ in direction $e$ if $l_F(e)\cap F=l_F(e)\cap P$. Denote by ${\mathcal F}_e(P)$ a set of all faces of $P$ that belong to the shadow boundary of $P$ in direction of $e$.

Note that the face $F$ is transformed into a face $F+z(e)$ in the extension $P^e=P+z(e)$. Denote dimension of $F$ by dim$F$. Lemma~\ref{F+z} below helps to understand how change faces of $P^e$ with respect to faces of $P$. Assertions of Lemma~\ref{F+z} are obvious.
\begin{lem}
\label{F+z}
Let $F$ be a face of a polytope $P$. Consider the sum $P^e=P+z(e)$.
There are the following three possibilities for the sum $F+z(e)$ :

(i) if $F$ is parallel to $e$, then $F+z(e)=F^e$ is an extension of $F$, and dim$(F+z(e))={\rm dim}F$;

(ii)  if $F$ is transversal to $e$ and $F\notin{\mathcal F}_e(P)$, then $F+z(e)$ is a parallel shift of $F$;

(iii)  if $F$ is transversal to $e$ and $F\in{\mathcal F}_e(P)$, then $F+z(e)$ is direct sum of $F$ and $z(e)$, and dim$(F+z(e))={\rm dim}F+1$. \hfill $\Box$
\end{lem}

According to Lemma~\ref{F+z}, each {\em facet} $F$ of the sum $P+bz(e)$ has one of the following three types

(i) extension $F=F_1^e$ of a facet $F_1$ of $P$;

(ii) a parallel shift $F=F_1+bz(e)$ of a facet $F_1$ of $P$;

(iii) direct sum $F=G\oplus bz(e)$ of a $(d-2)$-face $G$ of $P$ and the segment $bz(e)$.

Now consider Minkowski sum of a polytope $P(a)$ given in (\ref{Pa}) and a  segment $bz(e)$. Suppose that signs of $\langle p,e\rangle$ take all three values of the set $\{0,\pm 1\}$ for all $p\in\mathcal P$. Then, by Lemma~\ref{l2}, $bz(e)=P(f_e)$, where the function $f_e(p)$ is defined in (\ref{fe}). Obviously, the polytope $P(a)$ is supported by the hyperplane $H_p(a)$ defined in (\ref{Hpa}) at every {\em facet} of $P$.
\begin{prop}
\label{prs}
Let $P=P(a)$ be a polytope described in (\ref{Pa}). Suppose that each $(d-2)$-face $F\in{\mathcal F}_e(P)$, which is transversal to $e$, is a contact face $F=F(p)$ such that $\langle e,p\rangle=0$. Then the following equalities hold
\[P(a)+bz(e)=P(a)+P(f_e)=P(a+f_e).\]
\end{prop}
{\bf Proof}. We show that each facet $F_e(p)$ of the polytope $P^e=P(a)+bz(e)$ is supported by the hyperplane $H_p(a+f_e)$.

By Lemma~\ref{l1}, the hyperplane $H_p(f_e)$ supports the segment $bz(e)$ for any $p\in\mathcal P$. Similarly, $H_p(a)$ supports $P(a)$ if $p$ is a normal vector of a facet of $P(a)$. Let $F(p)$ be a contact face of $P(a)$. Consider the above three cases (i), (ii) and (iii) of facets of the sum $F(p)+bz(e)$.

Case (i). Let $F(p)\in{\mathcal F}_e(P)$ be a facet. Then $e\parallel F(p)$ and therefore $\langle p,e\rangle=0$. Hence $f_e(p)=0$, and $H_p(a+f_e)=H_p(a)$ supports the facet $F_e(p)= F(p)+bz(e)=F^e(p)$ of the sum $P+bz(e)$.

Case (ii). Let $F(p)$ be a facet and  $F(p)\notin{\mathcal F}_e(P)$. Then $\langle p,e\rangle\not=0$. Hence the sum $F(p)+bz(e)$ is a shift of $F(p)$ obtained as follows. Let $x\in F(p)$. Then the point
\[x+be\frac{\langle p,e\rangle}{|\langle p,e\rangle|} \]
belongs to $F(p)+bz(e)$. Here the multiple $\frac{\langle p,e\rangle} {|\langle p,e\rangle|}$ describes direction of the shift. Since $F(p)$ is a facet of $P(a)$, we have $\langle p,x\rangle=a(p)$ and therefore
\[\langle p,x+be\frac{\langle p,e\rangle}{|\langle p,e\rangle|}\rangle=a(p)+f_e(p)=(a+f_e)(p).  \]
Since $x$ is an arbitrary point of $F(p)$, this implies that the facet $F_e(p)=F(p)+bz(e)$ of $P+bz(e)$ is supported by $H_p(a+f_e)$.

Case (iii). Now, let $F(p)$ be a $(d-2)$-face of $P=P(a)$, $F(p)\in{\mathcal F}_e(P)$ and $F(p)$ is transversal to $e$. Then we have the case (iii) of Lemma~\ref{F+z}. The face $F(p)$ is transformed into the facet $F_e(p)=F(p)\oplus bz(e)$ of $P+bz(e)$. Since $F(p)\in{\mathcal F}_e(P)$, we have $\langle p,e\rangle=0$. Hence $f_e(p)=0$ and the hyperplane $H_p(a+f_e)=H_p(a)$ supports the facet $F_e(p)=F(p)\oplus bz(e)$.

So, each facet of the sum $P+bz(e)$ is supported by a hyperplane $H_p(a+f_e)$. Hence, $P(a)+bz(e)\supseteq P(a+f_e)$. Since, by Lemma~\ref{l2}, $bz(e)=P(f_e)$, according to Lemma~\ref{a12}, we obtain assertion of this Proposition. \hfill $\Box$

\section{Minkowski sum of a Voronoi parallelotope with a segment}
Now consider Minkowski sum of a {\em Voronoi parallelotope} $P(a)$ and the segment $bz(e)$. Without loss of generality, we can suppose that $\mathcal P$ contains the set of contact vectors of all contact faces of $P(a)$.

Obviously, each segment is a parallelotope, and moreover a Voronoi parallelotope. In fact, we can choose lengths of vectors $p\in\mathcal P$ such that $\langle p,e\rangle\in\{0,\pm 1\}$. In this case the function $f_e(p)$ transforms in the quadratic form $a_e(p)$ defined in (\ref{ae}), and $P(a_e)$ is a Voronoi parallelotope.

Let ${\mathcal P}_s(a)\subseteq{\mathcal P}(a)$ be a set of contact vectors of facets of $P(a)$. They are normal vectors of facets. If, for a vector $e\in{\mathbb R}^d$, the inclusions $\langle p,e\rangle\in\{0,\pm w\}$ hold for all $p\in{\mathcal P}_s(a)$, then one can change a value of $b$ and the length of $e$ such that the length of the segment $bz(e)$ does not change and $\langle p,e\rangle \in\{0,\pm 1\}$ for all $p\in{\mathcal P}_s(a)$. Hence we will consider vectors $e\in{\mathcal P}_s^*(a)$, where the dual ${\mathcal P}_s^*(a)$ is defined in (\ref{P*}). Of course, there may be another vectors $p\in\mathcal P$ with $\langle p,e\rangle \in\{0,\pm 1\}$. Hence we introduce the following set
\begin{equation}
\label{Pe}
{\mathcal P}_e=\{p\in{\mathcal P}:\langle p,e\rangle\in\{0,\pm 1\}\}.
\end{equation}
\begin{lem}
\label{lae}
Let $e\in{\mathcal P}^*_s(a)$, and scalar products $\langle p,e\rangle$ take all three values $0,+1,-1$ for $p\in{\mathcal P}_e$. Then
\[bz(e)=P(a_e), \]
where the quadratic form $a_e$ is defined in (\ref{ae}), and $P(a)$ is defined in (\ref{Pa}).
\end{lem}
{\bf Proof}. It is easy to see that $f_e(p)=a_e(p)$ for all $p\in{\mathcal P}_e$. By Lemma~\ref{l1}, the hyperplane $H_p(a_e)$ supports the segment $bz(e)$ for all $p\in{\mathcal P}_e$. By Lemma~\ref{l2}, $bz(e)=P(a_e)$. \hfill $\Box$

\begin{lem}
\label{l8}
Let ${\mathcal P}_s(a)$ be a set of normal vectors of a Voronoi parallelotope $P=P(a)$. Let $e\in{\mathcal P}^*_s(a)$, and let $F\in{\mathcal F}_e(P)$ be a $(d-2)$-face of $P$ that is transversal to $e$. Then $F=F(p)$ is a contact face for a contact vector $p$ such that $\langle p,e\rangle=0$ .
\end{lem}
{\bf Proof}. Suppose to the contrary that $F$ generates a 6-belt $B$. Let $\pm p_1,\pm p_2,\pm p_3\in{\mathcal P}_s(a)$ be normal vectors of the 6-belt $B$. Let $F=F(p_1)\cap F(p_2)$. Note that $F(p_1),F(p_2) \notin{\mathcal F}_e(P)$, since $F$ is transversal to $e$ and $F\in{\mathcal F}_e(P)$. Hence, for $i=1,2$, $\langle p_i,e\rangle\not=0$, and therefore $\langle p_i,e\rangle\in\{\pm 1\}$. Since $e\in{\mathcal P}^*_s(a)$, without loss of generality, we can suppose that $\langle p_1,e\rangle=1$ and $\langle p_2,e\rangle=-1$. Let $F(p_3)\not=F(p_2)$ be the second facet of the 6-belt $B$ that is adjacent to $F(p_1)$. Since $P=P(a)$ is a Voronoi parallelotope, the equality $p_3=p_1-p_2$ holds. This equality implies the equality $\langle p_3,e\rangle=\langle p_1,e\rangle-\langle p_2,e\rangle=2$ that contradicts to $\langle p_3,e\rangle\in\{0,\pm 1\}$. Hence $F$ cannot generate a 6-belt. Therefore $F=F(p)$ is a contact face.

Obviously, $F(p)=F(p_1)\cap F(p_2)$, where $F(p_1),F(p_2)$ are facets of $P(a)$. Then $p=p_1+p_2$. Since $e\in{\mathcal P}_s^*(a)$, $\langle p_1,e\rangle,\langle p_2,e\rangle\in\{0,\pm 1\}$. If $\langle p_1,e\rangle=\langle p_2,e\rangle=0$, then $\langle p,e\rangle=0$. Otherwise, without loss of generality, we can suppose that $\langle p_1,e\rangle=1$, $\langle p_2,e\rangle=-1$. this implies $\langle p,e\rangle=0$. \hfill $\Box$

\begin{prop}
\label{pr1}
Let $P=P(a)$ be a Voronoi parallelotope defined in (\ref{Pa}), where ${\mathcal P}\supseteq{\mathcal P}_e\supseteq{\mathcal P}_s(a)$. Then
\[P(a)+bz(e)=P(a)+P(a_e)=P(a+a_e). \]
\end{prop}
{\bf Proof}. Let $F(p)\in{\mathcal F}_e(P)$ be a $(d-2)$-face that is transversal to $e$. Then, by Lemma~\ref{l8}, $F=F(p)$ is a contact face of $P=P(a)$ and $\langle p,e\rangle=0$. Since $a_e(p)=f_e(p)$ for all $p\in{\mathcal P}_e$, we can apply Proposition~\ref{prs}. Hence the assertion of this Proposition holds. \hfill $\Box$

\vspace{2mm}
So, we have proved that $P(a)+P(a_e)=P(a+a_e)$ if $e\in{\mathcal P}^*_s(a)$. Recall that $a_e(p)=b\langle p,e\rangle^2$. Now we will prove that if $P(a)$ is irreducible and the sum $P(a)+P(a_e)$ is a parallelotope, then one can choose the length of $e$ such that $e\in{\mathcal P}^*_s(a)$.

Let $G(P)$ be a graph whose vertices correspond to facets of $P$. Let $v(F)$ be a vertex of $G(P)$ related to a facet $F$. Two vertices $v(F_1)$ and $v(F_2)$ are adjacent in $G(P)$ if and only if $F_1\cap F_2$ is a $(d-2)$-face that generates a 6-belt of $P$. It is proved in \cite{Or} that the graph $G(P)$ is connected if and only if the parallelotope $P$ is irreducible. In particular, this means that, for each $p\in{\mathcal P}_s(a)$, the facet $F(p)$ belongs to a 6-belt.

\begin{lem}
\label{wp}
Let $P=P(a)$ be an irreducible parallelotope. Let $e\in{\mathbb R}^d$ be a vector such the sum $P^e=P+bz(e)$ is a parallelotope. Then one can choose a length of the vector $e$ such that $e\in{\mathcal P}_s^*(a)$.
\end{lem}
{\bf Proof}. Let
\[{\mathcal P}_s^+(a)=\{p\in{\mathcal P}_s(a):\langle p,e\rangle>0\}. \]
We show that there is a number $w>0$ such that $\langle p,e\rangle=w$ for all vectors $p\in{\mathcal P}_s^+(a)$.

Suppose that there are two vectors $p_1,p_2\in{\mathcal P}_s^+(a)$ such that $\langle p_1,e\rangle\not=\langle p_2,e\rangle$. Since the graph $G(P)$ is connected, there are $p,p'\in{\mathcal P}_s^+(a)$ such that the facets $F(p)$ and $F(p')$ belong to the same 6-belt $B$, but $\langle p,e \rangle\not=\langle p',e\rangle$. Since $P=P(a)$ is a Voronoi parallelotope, the vector $p-p'$ is a normal vector of a facet of the 6-belt $B$. Hence the 6-belt $B$ consists of facets $F(\pm p),F(\pm p')$ and $F(\pm(p-p'))$. Since the sum $P+bz(e)$ is a parallelotope, one pair of opposite facets of the belt $B$ belongs to the shadow boundary ${\mathcal F}_e(P)$. Normal vectors of these facets are orthogonal to the vector $e$. Since $\langle p,e\rangle \not=0$ and $\langle p',e\rangle \not=0$, we have $\langle p-p',e\rangle =0$. This contradicts to the above assertion that $\langle p,e \rangle\not=\langle p',e\rangle$. Therefore, $\langle p,e\rangle=w>0$ for all $p\in{\mathcal P}_s^+(a)$. Hence, one can choose a length of $e$ such that $\langle p,e\rangle\in\{0,\pm 1\}$ for all $p\in{\mathcal P}_s(a)$. \hfill $\Box$

\vspace{2mm}
Now we can prove Theorem~\ref{thm}.
\begin{theor}
\label{thm}
Let $P(a)$ be an irreducible Voronoi parallelotope, defined in (\ref{Pa}), where ${\mathcal P}\supseteq{\mathcal P}(a)\supseteq{\mathcal P}_s(a)$. Let $e\in{\mathbb R}^d$ be a vector. Then one can choose a length of the vector $e$ such that the following assertions are equivalent:

(i) Minkowski sum $P(a)+bz(e)$ is a Voronoi parallelotope for any $b\ge 0$, and
\[P(a)+bz(e)=P(a)+P(a_e)=P(a+a_e);\]

(ii) $e\in{\mathcal P}_s^*(a)$;
\end{theor}
{\bf Proof}. Lemma~\ref{wp} and definition (\ref{P*}) of the set ${\mathcal P}_s^*(a)$ imply the implication (i)$\Rightarrow$(ii).

We prove implication (ii)$\Rightarrow$(i). Without loss of generality, we can suppose that $e\in{\mathcal P}_s^*(a)\subseteq{\mathcal P}_e$. By Proposition~\ref{pr1}, $P(a)+bz(e)=P(a+a_e)$ is a Voronoi parallelotope. \hfill $\Box$

\vspace{2mm}
Theorem~\ref{thm} is a generalization of results for Voronoi polytopes of root lattices $D_n$, $E_6$ and $E_7$ obtained in papers \cite{Gr2}, \cite{DGM} and \cite{Gr3}, respectively.

Note that if the Voronoi parallelotope is reducible, then one can apply Theorem~\ref{thm} to each component separately.

Theorem~\ref{thm} has the following important Corollary.
\begin{cor}
\label{cr1}
If ${\mathcal P}_s^*(a)=\emptyset$, then $P(a)+bz(e)$ is not a parallelotope for any vector $e$.
\end{cor}

Examples of $P(a)$ with ${\mathcal P}_s^*(a)=\emptyset$ are Voronoi parallelotopes of dual root lattices $E_6^*$ and $E_7^*$ (see \cite{DGM} and \cite{Gr3}).

\end{document}